\documentclass[a4paper,11pt]{article}

\usepackage{amsmath}
\usepackage{amsfonts}
\usepackage{amsthm}
\usepackage{amssymb}

\title{Divisorial contractions in dimension three\\
which contract divisors to compound $A_1$ points}
\author{Masayuki Kawakita}
\date{}
\theoremstyle{definition}
\newtheorem{Definition}{Definition}[section]
\newtheorem{Remark}[Definition]{Remark}
\newtheorem{Notation}[Definition]{Notation}
\newtheorem{Construction}[Definition]{Construction}

\theoremstyle{plain}
\newtheorem{Theorem}[Definition]{Theorem}
\newtheorem{Proposition}[Definition]{Proposition}
\newtheorem{Lemma}[Definition]{Lemma}
\newtheorem{Corollary}[Definition]{Corollary}
\newtheorem{Claim}[Definition]{Claim}

\newcommand{\cO}{{\mathcal O}}

\newcommand{\cHom}{{\mathcal H}om}

\newcommand{\bC}{{\mathbb C}}
\newcommand{\bP}{{\mathbb P}}
\newcommand{\bQ}{{\mathbb Q}}
\newcommand{\bZ}{{\mathbb Z}}

\newcommand{\fm}{{\mathfrak m}}

\newcommand{\sM}{{\mathsf M}}

\newcommand{\Spec}{\mathrm{Spec}}

\newcommand{\divis}{\mathrm{div}}
\newcommand{\Ker}{\mathrm{Ker}}
\newcommand{\ord}{\mathrm{ord}}
\newcommand{\wt}{\mathrm{wt}}
\newcommand{\Bl}{\mathrm{Bl}}
\newcommand{\LLC}{\mathrm{LLC}}

\newcommand{\others}{{(\mathrm{others})}}

\begin{document}
\maketitle

\begin{abstract}
We deal with a divisorial contraction in dimension $3$ which
contracts its exceptional divisor to a c$A_1$ point.
We prove that any such contraction is obtained by a suitable
weighted blow-up.
\end{abstract}

\section{Introduction}\label{sec:intro}
Explicit description of divisorial contractions
is a beautiful object in itself,
and in dimension $3$ it is one of the most important remaining problems.
The aim of this paper is to continue the study of this,
following my previous paper \cite{Kwk00}.

Let $f \colon (Y \supset E) \to (X \ni P)$ be
a divisorial contraction in dimension $3$
which contracts its exceptional divisor $E$ to a point $P$.
The theorem in \cite{Kwk00} is that any such contraction to a smooth point $P$
is obtained by a suitable weighted blow-up.
In the proof of this theorem, a numerical game for types of
singularities on $Y$ and for dimensions of $\cO_X/f_*\cO_Y(-iE)$'s
plays one of essential roles, and it works also even if
$P$ is a Gorenstein singularity.
In this paper we treat the case where $P$ is a c$A_1$ point,
starting with this game, and prove the following theorem.

\begin{Theorem}\label{th:main-th-intro}
\textup{($=$ \textbf{Theorem \ref{th:main-th}.})}
Let $Y$ be a $3$-dimensional $\bQ$-factorial normal
variety with only terminal singularities,
and let $f \colon (Y \supset E) \to (X \ni P)$ be an algebraic germ of
a divisorial contraction which contracts its exceptional divisor
$E$ to a c$A_1$ point $P$.
Then $f$ is obtained by a suitable weighted blow-up.
More precisely, under a suitable  analytic identification
$P \in X \cong o \in (xy+z^2+w^N=0) \subset \bC^4$,
$f$ is one of the following weighted blow-ups.

\smallskip
\noindent \textup{(\ref{th:main-th-intro}.1)}
$f$ is the weighted blow-up with its weights
$\wt(x,y,z,w)=(s,2t-s,t,1)$,
where $s,t$ are coprime positive integers such that $s \le t \le N/2$.\\
\textup{(\ref{th:main-th-intro}.2)}
$N=3$ and $f$ is the weighted blow-up
with its weights $\wt(x,y,z,w)=(1,5,3,2)$.
\end{Theorem}

The hardest part of this theorem lives in case (\ref{th:main-th-intro}.1).
Adding that there exist infinitely many such weighted blow-ups by the choice
of an analytic identification
$P \in X \cong o \in (xy+z^2+w^N=0) \subset \bC^4$,
some difficulties arise in controlling the value of $N$,
which should be large compared to the discrepancy of $f$.
For this, we introduce a special surface $P \in S \subset X$
(Definition \ref{def:sp-surf}) and reduce the problem to constructing
a special surface of which strict transform on $Y$ has 
only relatively mild singularities.

Y.\@ Kawamata has succeeded the description in the case where
$P$ is a terminal quotient singularity (\cite{Kwm96}),
and A.\@ Corti has done in the case
where $P$ is an ordinary double point (\cite[Theorem 3.10]{Co00}),
a special case of Theorem \ref{th:main-th-intro}.
In this paper we can also see the essence of their proof,
comparing discrepancies and using Shokurov's
connectedness lemma respectively.

I would like to thank Prof.\@ Yujiro Kawamata and Prof.\@ Alessio Corti
for their stimulating encouragement.
I am also grateful to Dr.\@ Nobuyuki Kakimi for his invaluable comments.
He told me the existence of weighted blow-ups in Theorem \ref{th:main-th-intro}
as examples of divisorial contractions.


\section{Statement of Theorem}\label{sec:state-th}
We work over the complex number field $\bC$.
A variety means an irreducible, reduced, separated scheme
of finite type over $\Spec \,\bC$.
Though our objects are algebraic in themselves
and we do in the algebraic category throughout the paper,
we often use analytic functions for convenience.
This produces no problem
by adding higher terms to them if necessary
to put them into algebraic functions.
Our argument depends not on the local ring $\cO_{X,P}$ itself,
but only on a quotient $\cO_{X,P}/\fm_P^n$ by a sufficiently large multiple
of the maximal ideal $\fm_P \subset \cO_{X,P}$.
We use basic terminologies in \cite[Chapters 1, 2]{K+92}.

First we define a divisorial contraction.
In this paper it means a morphism which may emerge in
the minimal model program.

\begin{Definition}\label{def:div-contr}
Let $f \colon Y \to X$ be a morphism with connected fibers between
normal varieties. We call $f$ a \textit{divisorial contraction}
if it satisfies the following conditions.

\smallskip
\noindent (\ref{def:div-contr}.1)
$Y$ is $\bQ$-factorial with only terminal singularities.\\
(\ref{def:div-contr}.2)
The exceptional locus of $f$ is a prime divisor.\\
(\ref{def:div-contr}.3)
$-K_Y$ is $f$-ample.\\
(\ref{def:div-contr}.4)
The relative Picard number of $f$ is $1$.
\end{Definition}

We recall the classification of terminal Gorenstein singularities in
dimension $3$.

\begin{Definition}\label{def:cDV-sing}
Let $P \in X$ be an algebraic germ (resp.\@ an analytic germ) of
a variety (resp.\@ an analytic space) in dimension $3$.
We call $P$ a \textit{cDV} (\textit{compound Du Val}) \textit{point}
if a general hyperplane section is normal and has a Du Val 
singularity at $P$.
The singularity $P$ is said to be
\textit{c$A_n$, c$D_n$, c$E_n$} (\textit{compound $A_n, D_n, E_n$})
according to the type of the Du Val singularity on
a general hyperplane section.
\end{Definition}

\begin{Theorem}\label{th:term-Gor}
\textup{(\cite[Theorem1.1]{R83}.)}
Let $P \in X$ be an algebraic germ
\textup{(}resp.\@ an analytic germ\textup{)} of
a  normal variety \textup{(}resp.\@ analytic space\textup{)}
in dimension $3$.
Then $P$ is a terminal Gorenstein singularity if and only if
$P$ is an isolated cDV point.
\end{Theorem}

\begin{Remark}\label{rem:th:term-Gor}
(\ref{rem:th:term-Gor}.1)
Let $P \in X \cong o \in (f=0) \subset \bC^4$ be
a terminal Gorenstein singularity in dimension $3$.
We can divide such singularities by the rank $r$
of the Hessian matrix of $f$ at $o$:

\smallskip
\noindent (\ref{rem:th:term-Gor}.1.1)
$r=1$. $P$ is c$D_n$, c$E_6$, c$E_7$, or c$E_8$.\\
(\ref{rem:th:term-Gor}.1.2)
$r=2$. $P$ is c$A_n$ with $n \ge 2$.\\
(\ref{rem:th:term-Gor}.1.3)
$r=3$. $P$ is c$A_1$, but is not an ordinary double point.\\
(\ref{rem:th:term-Gor}.1.4)
$r=4$. $P$ is an ordinary double point.

\smallskip
\noindent (\ref{rem:th:term-Gor}.2)
If $P$ is an isolated c$A_1$ point,
we have an analytic identification
$P \in X \cong o \in (xy+z^2+w^N=0) \subset \bC^4$ for some $N \ge 2$.
This $N$ depends only on $P \in X$ itself.
\end{Remark}

Now it is the time when we state the theorem precisely.

\begin{Theorem}\label{th:main-th}
Let $Y$ be a $3$-dimensional $\bQ$-factorial normal
variety with only terminal singularities,
and let $f \colon (Y \supset E) \to (X \ni P)$ be an algebraic germ of
a divisorial contraction which contracts its exceptional divisor
$E$ to a c$A_1$ point $P$.
Then $f$ is obtained by a suitable weighted blow-up.
More precisely, under a suitable  analytic identification
$P \in X \cong o \in (xy+z^2+w^N=0) \subset \bC^4$,
$f$ is one of the following weighted blow-ups.

\smallskip
\noindent \textup{(\ref{th:main-th}.1)}
$f$ is the weighted blow-up with its weights
$\wt(x,y,z,w)=(s,2t-s,t,1)$,
where $s,t$ are coprime positive integers such that $s \le t \le N/2$.\\
\textup{(\ref{th:main-th}.2)}
$N=3$ and $f$ is the weighted blow-up
with its weights $\wt(x,y,z,w)=(1,5,3,2)$.
\end{Theorem}

\begin{Remark}\label{rem:th:main-th}
Consider an analytic germ of a c$A_1$ point
$o \in (xy+z^2+w^N=0) \subset \bC^4$ ($N \ge 2$)
and blow-up this with weights as one of them in Theorem \ref{th:main-th}.
Then the exceptional locus of this weighted blow-up is irreducible, and
the weighted blown-up analytic space has actually
only terminal singularities.
\end{Remark}


\section{Singular Riemann-Roch Technique}\label{sec:singRR}
In this section we state some numerical results obtained by using the singular
Riemann-Roch formula (\cite[Theorem 10.2]{R87}),
most of which are in \cite{Kwk00}.
Let $Y$ be a $3$-dimensional $\bQ$-factorial normal
variety with only terminal singularities,
and let $f \colon (Y \supset E) \to (X \ni P)$ be an algebraic germ of
a divisorial contraction which contracts its exceptional divisor $E$
to a Gorenstein point $P$. 
Throughout this section we fix this situation and spread a general theory.

Let $K_Y=f^*K_X+aE$ and let $r$ be the global Gorenstein index of $Y$,
that is, the smallest positive integer such that $rK_Y$ is Cartier.
Because $a$ and $r$ are coprime
by \cite[Lemma 2.5]{Kwk00}, we can take an integer $e$ such that
$ae \equiv 1$ modulo $r$.

Let $I = \{Q: \mathrm{type} \ \frac{1}{r_Q}(1, -1, b_Q)\}$
be the set of fictitious singularities of $Y$, that is,
terminal quotient singularities obtained by
flat deformations of non-Gorenstein singularities of $Y$.
Then $(\cO_{Y_Q}(E_Q))_Q \cong (\cO_{Y_Q}(eK_{Y_Q}))_Q$,
where $(Y_Q, E_Q)$ is the deformed pair near $Q$ from $(Y, E)$.
We note that $b_Q$ is coprime to $r_Q$ and that
$e$ is also coprime to $r_Q$ because $r|(ae-1)$.
Hence $v_Q= \overline{eb_Q}$ is coprime to $r_Q$.
Here $\bar{\ }$ denotes the smallest residue modulo $r_Q$,
that is, $\overline{j} = j - \lfloor \frac{j}{r_Q} \rfloor r_Q$,
where $\lfloor \ \rfloor$ denotes the round down, that is,
$\lfloor j \rfloor = \max\{k \in \bZ | k \le j\}$.
Replacing $b_Q$ with $r_Q - b_Q$ if necessary,
we may assume that $v_Q \le r_Q / 2$.
With this description, $r = 1$ if $I$ is empty, and otherwise
$r$ is the lowest common multiple of $\{r_Q\}_{Q \in I}$.
We put $J=\{(r_Q, v_Q)\}_{Q \in I}$.

\begin{Proposition}\label{prop:fromKwk}
\textup{(\ref{prop:fromKwk}.1)}
$rE^3 \in \bZ_{>0}$.\\
\textup{(\ref{prop:fromKwk}.2)}
$\displaystyle
aE^3 = 2 - \sum_{Q\in I} \frac{v_Q(r_Q - v_Q)}{r_Q}$.\\
\textup{(\ref{prop:fromKwk}.3)}
For $1 \le i \le a$,
\begin{align*}
\dim_{\bC}\cO_X/f_*\cO_Y(-iE) = i^2 - \frac{1}{2}\sum_{Q \in I}
\min_{0 \le j < i}\{(1+j)jr_Q + i(i-1-2j)v_Q\}.
\end{align*}
\textup{(\ref{prop:fromKwk}.4)}
If $a \ge 2$,
then $\displaystyle \sum_{Q \in I}v_Q = 3 - \dim_{\bC}\fm_P/f_*\cO_Y(-2E)$,\\
where $\fm_P \subset \cO_X$ is the ideal sheaf of $P \in X$.
\end{Proposition}

\begin{proof}
We see them in \cite[Proposition 2.7]{Kwk00}.
$a \ge 2$ is used in the proof of (D) in \cite[Proposition 2.7]{Kwk00}.
\end{proof}

Now we classify $f$ from the numerical point of view.
\begin{Theorem}\label{th:fromKwk}
Exactly one of the following holds.

\smallskip
\noindent \textup{(\ref{th:fromKwk}.1)}
$a=1$.\\
\textup{(\ref{th:fromKwk}.2)}
$a \ge 2$, and\\
\textup{(\ref{th:fromKwk}.2.0)}
$\dim_{\bC}\fm_P/f_*\cO_Y(-2E)=0$.\\
In this case $J=\{(7, 3)\}$ or $\{(3, 1), (5, 2)\}$, and $a=2$.\\
\textup{(\ref{th:fromKwk}.2.1)}
$\dim_{\bC}\fm_P/f_*\cO_Y(-2E)=1$, and\\
\textup{(\ref{th:fromKwk}.2.1.1)} $J=\{(r, 2)\}$.
In this case $a=2$ or $4$.\\
\textup{(\ref{th:fromKwk}.2.1.2)} $J=\{(r_1, 1), (r_2, 1)\}$ $(r_1 \le r_2)$.\\
\textup{(\ref{th:fromKwk}.2.2)}
$\dim_{\bC}\fm_P/f_*\cO_Y(-2E)=2$.\\
In this case $J=\{(r, 1)\}$.\\
\textup{(\ref{th:fromKwk}.2.3)}
$\dim_{\bC}\fm_P/f_*\cO_Y(-2E)=3$.\\
In this case $P$ is a smooth point, $f$ is the usual blow-up along $P$,
and $a=2$.
\end{Theorem}

\begin{proof}
Most of the results come from \cite{Kwk00}.
In the case $a \ge2$, we classify $J$ by Proposition \ref{prop:fromKwk}.4
according to the value of $\dim_{\bC}\fm_P/f_*\cO_Y(-2E)$.
In case (\ref{th:fromKwk}.2.1.1), we have $a=2$ or $4$ from $a(rE^3)=4$ and
$rE^3 \in \bZ$, which are obtained by Propositions \ref{prop:fromKwk}.1-2.
In case (\ref{th:fromKwk}.2.3), $Y$ is Gorenstein
and thus \cite[Theorem 5]{Cu88} induces the result.
In case (\ref{th:fromKwk}.2.0), we know
all the possible values of $J$ and the corresponding values of $aE^3$
as follows by \cite[Subsection 2.3]{Kwk00}.
\begin{center}
\begin{tabular}{|c|c|c|c|c|c|}
\hline
\multicolumn{2}{|c|}{$J=\{(r, 3)\}$} &
\multicolumn{2}{c|}{$J=\{(r_1, 1), (r_2, 2)\}$} &
\multicolumn{2}{c|}{$J=\{(r_1, 1), (r_2, 1), (r_3, 1)\}$} \\
\hline
$r$ & $aE^3$ & $(r_1, r_2)$ & $aE^3$ & $(r_1, r_2, r_3)$ & $aE^3$ \\
\hline
$7$ & $2/7$  & $(2, 5)$     & $3/10$ & $(2, 2, r_3)$     & $2/2r_3$ \\
$8$ & $1/8$  & $(3, 5)$     & $2/15$ & $(2, 3, 3)$       & $1/6$ \\
    &        & $(4, 5)$     & $1/20$ & $(2, 3, 4)$       & $1/12$ \\
    &        & $(2, 7)$     & $1/14$ & $(2, 3, 5)$       & $1/30$ \\
\hline
\end{tabular}
\end{center}
Considering that $a \ge 2$, $rE^3 \in \bZ$,
and that $a, r_Q$ are coprime,
we can restrict the possibility again to the three cases below.

\smallskip
\begin{tabular}{l}
$J=\{(7,3)\}$ and $a=2$.\\
$J=\{(2,1), (5,2)\}$ and $a=3$.\\
$J=\{(3,1), (5,2)\}$ and $a=2$.\\
\end{tabular}

\smallskip
\noindent Thus we have only to exclude the case $J=\{(2,1), (5,2)\}$.

In the case $J=\{(2,1), (5,2)\}$, we have $a=3$,
$e \equiv 7$ modulo $10$, $E^3=1/10$, and
$\{(r_Q,v_Q,b_Q)\}_{Q \in I} = \{(2,1,1), (5,2,1)\}$.
From the proof of \cite[Proposition 2.7]{Kwk00} we know that, for $i \le a$,
\begin{align*}
&\dim_{\bC}f_*\cO_Y(iE)/f_*\cO_Y((i-1)E)\\
&= \frac{1}{12} \{2(3i^2-3i+1)
- 3(2i-1)a + a^2\}E^3 + \frac{1}{12} E\cdot c_2(Y) + A_i - A_{i-1},
\end{align*}
where $\displaystyle
A_i = \sum_{Q\in I}
\Bigl( - \overline{ie} \frac{r_Q^2-1}{12r_Q} +
\sum_{j=1}^{\overline{ie}-1}
\frac{\overline{jb_Q}(r_Q-\overline{jb_Q})}{2r_Q} \Bigr)$.\\
In our case we have
\begin{align*}
&\dim_{\bC}f_*\cO_Y(iE)/f_*\cO_Y((i-1)E) \\
&= \frac{1}{60}(3i^2-12i+10) + \frac{1}{12} E\cdot c_2(Y) + A_i - A_{i-1}
\quad (i \le 3),
\end{align*}
$A_1 = -21/40$, $A_2 = 0$, and so on.
Putting $i=1,2$ we obtain
\begin{alignat*}{2}
0 &= \frac{1}{60} + \frac{1}{12} E\cdot c_2(Y) - \frac{21}{40}
&\quad&(i=1),\\
0 &= - \frac{1}{30} + \frac{1}{12} E\cdot c_2(Y) + \frac{21}{40}
&&(i=2).
\end{alignat*}
These two equations contradict each other.
\end{proof}


\section{First Step to Proof}\label{sec:first}
In this section we take the first step to the proof of
Theorem \ref{th:main-th}.
We keep numerical data in Section \ref{sec:singRR}.

\smallskip
\begin{tabular}{l}
$K_Y=f^*K_X+aE$ \\
$r$ : the Gorenstein index \\
$e$ : an integer such that $ae \equiv 1$ modulo $r$ \\
$I = \{Q: \mathrm{type} \ \frac{1}{r_Q}(1, -1, b_Q)\}$
: the set of fictitious singularities of $Y$ \\
$v_Q = \overline{eb_Q}$ \ $(v_Q \le r_Q/2)$ \\
$J=\{(r_Q, v_Q)\}_{Q \in I}$
\end{tabular}

\smallskip
\noindent Additionally, we define an integer $N \ge 2$ as follows.
\begin{align*}
P \in X \cong o \in (xy+z^2+w^N=0) \subset \bC^4.
\end{align*}

First we construct a tower of normal varieties.

\begin{Construction}\label{constr:tower}
We construct birational morphisms $g_i \colon  X_i \to X_{i-1}$
between normal factorial varieties,
irreducible and reduced closed subschemes $Z_i \subset X_i$,
and prime divisors $F_i$ on $X_i$ inductively,
and define positive integers $n, m$, with the following procedure.

\smallskip
\noindent (\ref{constr:tower}.1)
Define $X_0$ as $X$ and $Z_0$ as $P$.\\
(\ref{constr:tower}.2.1)
If $Z_{i-1}$ is a point,
we define $g_i$ as the blow-up of $X_{i-1}$ along $Z_{i-1}$.\\
(\ref{constr:tower}.2.2)
If $Z_{i-1}$ is a curve,
we define $b_i \colon \Bl_{Z_{i-1}}(X_{i-1}) \to X_{i-1}$ as
the blow-up of $X_{i-1}$ along $Z_{i-1}$, and define
${b'}_{\!i} \colon X_i \to \Bl_{Z_{i-1}}(X_{i-1})$ as a resolution of
singularities near $b_i^{-1}(Z_{i-1})$.
Precisely, ${b'}_{\!i}$ is a proper morphism which is isomorphic
over $\Bl_{Z_{i-1}}(X_{i-1}) \setminus b_i^{-1}(Z_{i-1})$, and
$X_i$ is smooth near $(b_i\circ {b'}_{\!i})^{-1}(Z_{i-1})$.
We note that ${b'}_{\!i}$ is isomorphic at the generic point of
the center of $E$ on $\Bl_{Z_{i-1}}(X_{i-1})$.
We define $g_i = b_i \circ {b'}_{\!i} \colon X_i \to X_{i-1}$.\\
(\ref{constr:tower}.3)
Define $Z_i$ as the center of $E$ on $X_i$ with the
reduced induced closed subscheme structure, and $F_i$ as
the only $g_i$-exceptional prime divisor on $X_i$ which contains $Z_i$.\\
(\ref{constr:tower}.4)
We stop this process when $Z_n = F_n$.
This process must terminate after finite steps like
\cite[Construction 2.1]{Kwk00} and thus we get the sequence
$X_n \to \cdots \to X_0$.\\
(\ref{constr:tower}.5)
We define $m \le n$ as the largest integer such that $Z_{m-1}$ is a point.\\
(\ref{constr:tower}.6)
We define $g_{ji}$ ($j \le i$) as the induced morphism from $X_i$ to $X_j$.
\end{Construction}

\begin{Remark}\label{rem:constr:tower-class}
$Z_i \subseteq F_i$ ($1 \le i\le n$) is exactly one of the following.

\smallskip
\noindent (\ref{rem:constr:tower-class}.1) $1 \le i < m$.
$Z_i$ is a point, and\\
(\ref{rem:constr:tower-class}.1.1)
$Z_i \in F_i \cong$ the vertex point $\in Q_0$,
where $Q_0$ denotes the cone $(xy+z^2 = 0) \subset \bP^3$
with homogeneous coordinates $x,y,z,w$.\\
(\ref{rem:constr:tower-class}.1.2)
$Z_i \in F_i \cong$ a non-vertex point $\in Q_0$.\\
(\ref{rem:constr:tower-class}.1.3)
$Z_i \in F_i \cong$ a point $\in Q$,
where $Q$ denotes the smooth quadratic $(xy+zw = 0) \subset \bP^3$
with homogeneous coordinates $x,y,z,w$.\\
(\ref{rem:constr:tower-class}.1.4)
$Z_i \in F_i \cong$ a point $\in \bP^2$.\\
(\ref{rem:constr:tower-class}.2) $m \le i < n$. $Z_i$ is a curve, and\\
(\ref{rem:constr:tower-class}.2.1) $i = m < n$ and,\\
(\ref{rem:constr:tower-class}.2.1.1)
$Z_i \subset F_i \cong$ a curve $\subset Q_0$.\\
(\ref{rem:constr:tower-class}.2.1.2)
$Z_i \subset F_i \cong$ a curve $\subset Q$.\\
(\ref{rem:constr:tower-class}.2.1.3)
$Z_i \subset F_i \cong$ a curve $\subset \bP^2$.\\
(\ref{rem:constr:tower-class}.2.2) $m < i < n$.\\
(\ref{rem:constr:tower-class}.3) $i=n$. $Z_i = F_i$ is a surface.
\end{Remark}

\begin{Remark}\label{rem:constr:tower-ideal}
We remark that $f_*\cO_Y(-iE) = g_{0n*}\cO_{X_n}(-iF_n)$ for any $i$
because $E$ and $F_n$ are the same as valuations.
\end{Remark}

By the next lemma, we have only to prove that
$F_n$ equals, as valuations, the only exceptional divisor obtained by
a weighted blow-up of $X$ emerging in Theorem \ref{th:main-th}.
This lemma is a generalization of \cite[Lemma 2.2]{Kwk00}, which
should be replaced by this.

\begin{Lemma}\label{lem:excep}
Let $f_i \colon Y_i \to X$ with $i = 1,2$ be
projective birational morphisms between normal varieties with
$f_i$-exceptional and $f_i$-ample $\bQ$-Cartier divisors $-E_i$ on $Y_i$.
Assume that $E_1$ and $E_2$ are the same as $\bQ$-linear combinations of
valuations.
Then $f_1$ and $f_2$ are isomorphic as morphisms over $X$.
\end{Lemma}

Second we evaluate various discrepancies and multiplicities.

\begin{Notation}\label{not:degree}
(\ref{not:degree}.1)
We define a positive integer $l \le m$ as the largest integer satisfying
that $l=1$ or that $Z_{l-1}$ is of type (\ref{rem:constr:tower-class}.1.1).\\
(\ref{not:degree}.2)
For curves $Z_i$ ($m \le i <n$),
we define the degree $d_i$ of $Z_i$ as follows.\\
(\ref{not:degree}.2.1)
In case (\ref{rem:constr:tower-class}.2.1.1-2),
$d_i$ denotes the degree of $Z_i$ considered as a subvariety in $\bP^3$
as in Remark \ref{rem:constr:tower-class}.\\
(\ref{not:degree}.2.2)
In case (\ref{rem:constr:tower-class}.2.1.3),
$d_i$ denotes the degree of $Z_i$ considered as a subvariety in $\bP^2$
as in Remark \ref{rem:constr:tower-class}.\\
(\ref{not:degree}.2.3)
In case (\ref{rem:constr:tower-class}.2.2),
$d_i$ denotes the degree of the finite morphism $Z_i \to Z_{i-1}$.
\end{Notation}

\begin{Notation}\label{not:first-discrep-mult}
Let $X$ be a normal variety and let $E$ be an algebraic
valuation, that is, a valuation
of the function field of $X$ which is obtained as an exceptional
divisor of some birational morphism $f \colon Y \to X$
from a normal variety $Y$.
Let $D$ be a $\bQ$-divisor on $X$ or $D = \mu \sM$ where $\mu$ is
a rational number and $\sM$ is a linear system
of finite dimension on $X$ which has no base points in codimension $1$.

\smallskip
\noindent (\ref{not:first-discrep-mult}.0)
Let $Z$ be a normal variety which is birational to $X$.
$D_Z$ denotes the strict transform of $D$ on $Z$.\\
(\ref{not:first-discrep-mult}.1)
Assume that $K_X+D$ is $\bQ$-Cartier.
$\alpha_{K_X + D}(E)$ denotes the discrepancy of $E$
with respect to $K_X + D$,
that is, $K_Y + D_Y = f^*(K_X + D) + \alpha_{K_X + D}(E)E
+ \others$.\\
(\ref{not:first-discrep-mult}.2)
Assume that $D$ is $\bQ$-Cartier.
$m_D(E)$ denotes the multiplicity of $E$ with respect to $D$,
that is, $f^*D = D_Y + m_D(E)E + \others$.
\end{Notation}

\begin{Notation}\label{not:discrep-mult}
Let $\sM$ be a general $f$-very ample linear system
of finite dimension on $Y$.
We define positive rational numbers $\mu, c_i$
by the following equations.

\smallskip
\noindent (\ref{not:discrep-mult}.1)
$\displaystyle
K_Y + \mu \sM = f^*(K_X + \mu \sM_X)$.\\
(\ref{not:discrep-mult}.2)
$\displaystyle
g_{0n}^*(\mu \sM_X) = \mu \sM_{X_n}
+ \sum_{1 \le i \le n} c_i(g_{in}^*F_i) + \others$.
\end{Notation}

\begin{Remark}\label{rem:not:discrep-mult}
(\ref{rem:not:discrep-mult}.1)
Because $\sM$ be a general $f$-very ample linear system on $Y$,
for any algebraic valuation $G$ we have
$\alpha_{K_X + \mu \sM_X}(G) = \alpha_{K_Y}(G)$.\\
(\ref{rem:not:discrep-mult}.2)
Putting $G = F_i$ in (\ref{rem:not:discrep-mult}.1), we obtain
\begin{align*}
\alpha_{K_X}(F_i) - \sum_{1 \le j \le i} c_jm_{F_j}(F_i) = \alpha_{K_Y}(F_i)
\begin{cases}
> 0 &(i < n), \\
= 0 &(i = n),
\end{cases}
\end{align*}
since $Y$ has only terminal singularities.
\end{Remark}

We give evaluation for $c_i$'s.

\begin{Proposition}\label{prop:eval}
\textup{(\ref{prop:eval}.1)}
$1 > c_1$ except the case $n=1$.\\
\textup{(\ref{prop:eval}.2)}
$c_n > \alpha_{K_{X_{n-1}}}(F_n)$ except the case $n=1$.\\
\textup{(\ref{prop:eval}.3.1)}
If $Z_i$ is a point of type
\textup{(\ref{rem:constr:tower-class}.1.1)}
or \textup{(\ref{rem:constr:tower-class}.1.4)},
then $c_i \ge c_{i+1}$.\\
\textup{(\ref{prop:eval}.3.2)}
If $Z_i$ is a point of type \textup{(\ref{rem:constr:tower-class}.1.2)},
then $2c_i \ge c_{i+1}$.\\
\textup{(\ref{prop:eval}.3.3)}
If $Z_i$ is a curve of type
\textup{(\ref{rem:constr:tower-class}.2.1.3)}
 or \textup{(\ref{rem:constr:tower-class}.2.2)},
then $c_i \ge d_ic_{i+1}$.\\
\textup{(\ref{prop:eval}.3.4)}
If $Z_i$ is a curve of type \textup{(\ref{rem:constr:tower-class}.2.1.1)},
then $2c_i \ge d_ic_{i+1}$.\\
\textup{(\ref{prop:eval}.4)}
If $Z_i$ is of type
\textup{(\ref{rem:constr:tower-class}.1.3)}
or \textup{(\ref{rem:constr:tower-class}.2.1.2)},
then $c_i \ge 1$.
\end{Proposition}

\begin{proof}
(\ref{prop:eval}.1)
Putting $i=1$ into (\ref{rem:not:discrep-mult}.2), we have $1-c_1>0$.

\smallskip
\noindent (\ref{prop:eval}.2)
We use Remark \ref{rem:not:discrep-mult}. Because
\begin{align*}
&\quad K_{X_n} + \mu\sM_{X_n} \\
&= g_n^*(K_{X_{n-1}} + \mu\sM_{X_{n-1}})
+ (\alpha_{K_{X_{n-1}}}(F_n) - c_n)F_n + \others \\
&=g_n^*(g_{0,n-1}^*(K_X + \mu\sM_X) + \alpha_{K_Y}(F_{n-1})F_{n-1}
+ \others) \\
&\quad + (\alpha_{K_{X_{n-1}}}(F_n) - c_n)F_n + \others \\
&=g_{0n}^*(K_X + \mu\sM_X) \\
&\quad +
(\alpha_{K_{X_{n-1}}}(F_n) - c_n + \alpha_{K_Y}(F_{n-1})m_{F_{n-1}}(F_n))F_n
+ \others,
\end{align*}
we have
$\alpha_{K_{X_{n-1}}}(F_n) - c_n + \alpha_{K_Y}(F_{n-1})m_{F_{n-1}}(F_n)
= \alpha_{K_Y}(F_n) = 0$.
Hence $c_n - \alpha_{K_{X_{n-1}}}(F_n) =
\alpha_{K_Y}(F_{n-1})m_{F_{n-1}}(F_n) > 0$.

\smallskip
\noindent (\ref{prop:eval}.3.1)
We will prove (\ref{prop:eval}.3) with the same idea.
Let $Z_i \in l \subset F_i$ be a general line
on $F_i \cong Q_0 \subset \bP^3$ or $\bP^2$ through $Z_i$,
and let $l'$ be its strict transform on $X_{i+1}$. Then,
\begin{align*}
0 \le (\mu\sM_{X_{i+1}} \cdot l')_{X_{i+1}}
= -c_{i+1}(F_{i+1} \cdot l')_{X_{i+1}} -c_i(F_i \cdot l)_{X_i}
= -c_{i+1} +c_i.
\end{align*}

\smallskip
\noindent (\ref{prop:eval}.3.2)
Let $Z_i \in c \subset F_i$ be a general conic
on $F_i \cong Q_0 \subset \bP^3$ through $Z_i$,
and let $c'$ be its strict transform on $X_{i+1}$. Then,
\begin{align*}
0 \le (\mu\sM_{X_{i+1}} \cdot c')_{X_{i+1}}
= -c_{i+1}(F_{i+1} \cdot c')_{X_{i+1}} -c_i(F_i \cdot c)_{X_i}
= -c_{i+1} +2c_i.
\end{align*}

\smallskip
\noindent (\ref{prop:eval}.3.3)
Let $l \subset F_i$ be a general line on $F_i \cong \bP^2$
in case (\ref{rem:constr:tower-class}.2.1.3)
and be a general fiber of $F_i \to Z_{i-1}$
in case (\ref{rem:constr:tower-class}.2.2),
and let $l'$ be its strict transform on $X_{i+1}$. Then,
\begin{align*}
0 \le (\mu\sM_{X_{i+1}} \cdot l')_{X_{i+1}}
= -c_{i+1}(F_{i+1} \cdot l')_{X_{i+1}} -c_i(F_i \cdot l)_{X_i}
= -d_ic_{i+1} +c_i.
\end{align*}

\smallskip
\noindent (\ref{prop:eval}.3.4)
Let $c \subset F_i$ be a general conic on $F_i \cong Q_0 \subset \bP^3$,
and let $c'$ be its strict transform on $X_{i+1}$. Then,
\begin{align*}
0 \le (\mu\sM_{X_{i+1}} \cdot c')_{X_{i+1}}
= -c_{i+1}(F_{i+1} \cdot c')_{X_{i+1}} -c_i(F_i \cdot c)_{X_i}
= -d_ic_{i+1} +2c_i.
\end{align*}

\smallskip
\noindent (\ref{prop:eval}.4)
The following proof is a generalization
of the proof of \cite[Theorem 3.10]{Co00} using
Shokurov's connectedness lemma (\cite[Theorem 17.4]{K+92}).
Let $P \in H \subset X$ be a general hyperplane section on $X$ through $P$,
and let $Z_{i-1} \in L \subset X_{i-1}$ be
a general hyperplane section on $X_{i-1}$ through $Z_{i-1}$
such that $Z_i \not\subseteq L_{X_i} \cap F_i$, and that
$L_{X_i} \cap F_i$ consists of
two lines $l_1 + l_2$ on $F_i \cong Q \subset \bP^3$,
which are fibers of two rulings of $Q \cong \bP^1 \times \bP^1$. Then
\begin{multline*}
g_{0i}^*(K_X+\mu \sM_X+\alpha_{K_Y}(F_{i-1})H)+g_i^*L \\
=K_{X_i}+\mu \sM_{X_i}+L_{X_i}+0{F_{i-1}}_{X_i}+c_iF_i+\others,
\end{multline*}
where we omit the term $\alpha_{K_Y}(F_{i-1})H$ if $i=1$.
Because
\begin{align*}
&\alpha_{(g_{0i}^*(K_X+\mu \sM_X+\alpha_{K_Y}(F_{i-1})H)+g_i^*L)}(F_n)\\
=&-m_{(\alpha_{K_Y}(F_{i-1})g_{0i}^*H+g_i^*L)}(F_n)
\le -m_L(F_n) = -1,
\end{align*}
we have
\begin{align*}
Z_i \subseteq
\LLC(X_i,\,g_{0i}^*(K_X+\mu \sM_X+\alpha_{K_Y}(F_{i-1})H)+g_i^*L),
\end{align*}
where LLC denotes the locus of log canonical singularities for a log pair,
that is,
the union of centers of all algebraic valuations with discrepancies
$\le -1$.
Moreover,
\begin{align*}
L_{X_i} \cap F_i \subseteq
\LLC(X_i,\,g_{0i}^*(K_X+\mu \sM_X+\alpha_{K_Y}(F_{i-1})H)+g_i^*L).
\end{align*}
Since $Z_i \not\subseteq L_{X_i} \cap F_i \cong l_1 + l_2$,
using the connectedness lemma for two small contractions 
contracting $l_1, l_2$ respectively, we obtain
\begin{align*}
\bP^1 \times \bP^1 \cong F_i \subseteq
\LLC(X_i,\,g_{0i}^*(K_X+\mu \sM_X+\alpha_{K_Y}(F_{i-1})H)+g_i^*L),
\end{align*}
that is, $c_i \ge 1$.
\end{proof}

We have a refined restriction as a corollary of preceding results.

\begin{Corollary}\label{cor:restr}
\textup{(\ref{cor:restr}.1)}
If $a=1$, then $f$ is the usual blow-up of $X$ along $P$.\\
\textup{(\ref{cor:restr}.2)}Assume that $a \ge 2$, that is, $n \ge 2$.\\
\textup{(\ref{cor:restr}.2.1.1)}
Case \textup{(\ref{th:fromKwk}.2.0)} never occurs.\\
\textup{(\ref{cor:restr}.2.1.2)}
Neither case \textup{(\ref{rem:constr:tower-class}.1.3)} nor
case \textup{(\ref{rem:constr:tower-class}.2.1.2)} occurs.\\
\textup{(\ref{cor:restr}.2.2)}
Exactly one of cases \textup{(\ref{rem:constr:tower-class}.1.2)} and
\textup{(\ref{rem:constr:tower-class}.2.1.1)} occurs.\\
\textup{(\ref{cor:restr}.2.3.1)}
$m<n$.\\
\textup{(\ref{cor:restr}.2.3.2)}
$\forall d_i=1$.\\
\textup{(\ref{cor:restr}.2.3.3)}
$2 > 2c_1 \ge \cdots \ge 2c_l \ge c_{l+1} \ge \cdots \ge c_n >1$.
\end{Corollary}

\begin{proof}
(\ref{cor:restr}.1)
This comes from Lemma \ref{lem:excep}.

\smallskip
\noindent(\ref{cor:restr}.2.1.1)
Since $a=2$ in case (\ref{th:fromKwk}.2.0), we have $n=2$ and

\smallskip
\begin{tabular}{l}
$Z_1$ is a point of type (\ref{rem:constr:tower-class}.1.1) and $N \ge 4$, or\\
$Z_1$ is a curve.
\end{tabular}

\smallskip
\noindent In both cases a general hyperplane section on $X$
through $P$ has multiplicity $1$ along $F_2$, which means that
$f_*\cO_Y(-2E)=g_{02*}\cO_{X_2}(-2F_2) \varsubsetneq \fm_P$.
This is a contradiction.

\smallskip
\noindent (\ref{cor:restr}.2.1.2)
Propositions \ref{prop:eval}.1, \ref{prop:eval}.3.1, and \ref{prop:eval}.4
imply this.

\smallskip
\noindent (\ref{cor:restr}.2.2)
If neither case (\ref{rem:constr:tower-class}.1.2) nor
case (\ref{rem:constr:tower-class}.2.1.1) occurs,
then from Proposition \ref{prop:eval} we have
\begin{align*}
1 > c_1 \ge (\prod_{m \le i <n}\!\!d_i)c_n >
(\prod_{m \le i <n}\!\!d_i)\alpha_{K_{X_{n-1}}}(F_n).
\end{align*}
This is a contradiction.

\smallskip
\noindent (\ref{cor:restr}.2.3.1-3)
We obtain them considering the following inequalities as
in the proof of (\ref{cor:restr}.2.2).
\begin{align*}
2 > 2c_1 \ge (\prod_{m \le i <n}\!\!d_i)c_n >
(\prod_{m \le i <n}\!\!d_i)\alpha_{K_{X_{n-1}}}(F_n).
\end{align*}
$m<n$ comes from $\alpha_{K_{X_{n-1}}}(F_n)=1$ and (\ref{cor:restr}.2.2).
\end{proof}

\begin{Remark}\label{rem:cor:restr}
Because Corollaries \ref{cor:restr}.2.2 and \ref{cor:restr}.2.3.2,
$F_l \cong Q_0$, $N \ge 2l+1$.
We define $Z_l \subseteq l_0 \subset F_l$ as the unique line
on $F_l \cong Q_0 \subset \bP^3$ containing $Z_l$.
\end{Remark}

The problem is reduced to investigating cases (\ref{th:fromKwk}.2.1-2),
which will be done in the following sections.
As the last part of this section, we give some information
for these remaining cases (\ref{th:fromKwk}.2.1-2).

\begin{Corollary}\label{cor:last}
\textup{(\ref{cor:last}.1)}
$Z_{i+1} \not\subseteq {F_i}_{X_{i+1}} \cap F_{i+1}$.\\
\textup{(\ref{cor:last}.2)}
$a=n+m-l$.\\
\textup{(\ref{cor:last}.3.1)}
In case \textup{(\ref{th:fromKwk}.2.1)},
$Z_1 \subset F_1 \cong Q_0$ in $\bP^3$ and it is a point.\\
\textup{(\ref{cor:last}.3.2)}
In case \textup{(\ref{th:fromKwk}.2.2)},
$Z_1 \subset F_1 \cong Q_0$ in $\bP^3$ and it is a line.\\
\end{Corollary}

\begin{proof}
(\ref{cor:last}.1)
This is trivial since $\fm_P \neq f_*\cO_Y(-2E) = g_{0n*}\cO_{X_n}(-2F_n)$.

\smallskip
\noindent (\ref{cor:last}.2)
This comes from Corollary \ref{cor:restr}.2.2 and (\ref{cor:last}.1).

\smallskip
\noindent (\ref{cor:last}.3)
$F_1 \cong Q_0$ comes from $a \ge 2$ and Corollary \ref{cor:restr}.2.1.2.
We know the shape of $Z_1 \subset F_1 \cong Q_0 \in \bP^3$
from the equation below.
\begin{align*}
&4 - \dim_{\bC}\fm_P/f_*\cO_Y(-2E) \\
=&\dim_{\bC}f_*\cO_Y(-2E)/\fm_P^2 \\
=&\dim_{\bC} \mathrm{Im}\,[(v \in \fm_P | Z_1 \subseteq \divis(v)_{X_1})
\to \fm_P/\fm_P^2] \\
=&\dim_{\bC}\{v \in \varGamma(\bP^3, \cO_{\bP^3}(1)) |
v = 0 \ \mathrm{or} \ Z_1 \subseteq \divis(v)\},
\end{align*}
where the second equality comes from $\fm_P \neq f_*\cO_Y(-2E)$.
\end{proof}


\section{Exceptional Case}\label{sec:excep}
In this section we treat the exceptional case (\ref{th:main-th}.2),
which corresponds to case (\ref{th:fromKwk}.2.1.1),
and our aim is the following.

\begin{Proposition}\label{prop:main-excep}
Assume that $f$ is of type \textup{(\ref{th:fromKwk}.2.1.1)}. Then
$f$ is obtained by a suitable weighted blow-up of type
\textup{(\ref{th:main-th}.2)}.
\end{Proposition}

Throughout this section, we assume that $f$ is of type (\ref{th:fromKwk}.2.1.1)
and struggle with Proposition \ref{prop:main-excep}.
We note that $1<m<n$ by the assumption and Corollaries \ref{cor:restr}.2.3.1
and \ref{cor:last}.3.1, and that $N \ge 3$ by Remark \ref{rem:cor:restr}.

First we restate the conclusion.

\begin{Lemma}\label{lem:restate-excep}
The following imply Proposition \textup{\ref{prop:main-excep}}.

\smallskip
\noindent \textup{(\ref{lem:restate-excep}.1)} $(n,m,l)=(3,2,1)$.\\
\textup{(\ref{lem:restate-excep}.2)}
$Z_2$ is a curve which intersects the strict transform of $l_0$ on $X_2$.
\end{Lemma}

\begin{proof}
Though analytic functions seem to emerge in this proof,
we stay in the algebraic category by adding higher terms to them if necessary,
as we have said in the first paragraph in Section \ref{sec:state-th}.
First we prove a claim on an analytic description.

\begin{Claim}\label{cl:lem:restate-excep-id}
There exists an identification
$P \in X \cong o \in (xy+z^2+w^N=0) \subset \bC^4$
satisfying the following conditions.

\smallskip
\noindent \textup{(\ref{cl:lem:restate-excep-id}.1)}
$l_0 = F_1 \cap \divis(y)_{X_1} \cap \divis(z)_{X_1}$.\\
\textup{(\ref{cl:lem:restate-excep-id}.2)}
$Z_1 = l_0 \cap \divis(w)_{X_1}$.\\
\textup{(\ref{cl:lem:restate-excep-id}.3)}
$Z_2 = F_2 \cap \divis(z)_{X_2}$.
\end{Claim}

\begin{proof}[Proof of Claim \textup{\ref{cl:lem:restate-excep-id}}]
It is trivial that we can choose an identification satisfying
(\ref{cl:lem:restate-excep-id}.1).
Then by (\ref{cl:lem:restate-excep-id}.1),
$Z_1 = l_0 \cap \divis(w+tx)_{X_1}$ for some $t \in \bC$.
Because $xy+z^2+w^N = xy'+z^2+(w')^N$ for
$w'=w+tx$ and $y'=y+(w^N-(w+tx)^N)/x$,
by replacing $y,w$ with $y',w'$
we may assume (\ref{cl:lem:restate-excep-id}.2) moreover.
Then $Z_2 = F_2 \cap \divis(z+tx^2)_{X_2}$ for some $t \in \bC$
by (\ref{lem:restate-excep}.2) and
Corollaries \ref{cor:restr}.2.3.2 and \ref{cor:last}.1. 
Because $xy+z^2+w^N = xy'+(z')^2+w^N$ for
$z'=z+tx^2$ and $y'=y-2txz-t^2x^3$,
by replacing $y,z$ with $y',z'$
we may assume (\ref{cl:lem:restate-excep-id}.3) moreover.
\end{proof}

Second we prove that $F_3$ equals, as valuations,
an exceptional divisor obtained by a weighted blow-up of $X$.

\begin{Claim}\label{cl:lem:restate-excep-wbl}
Under the identification in Claim \textup{\ref{cl:lem:restate-excep-id}},
$F_3$ equals, as valuations,
an exceptional divisor obtained by the weighted blow-up of $X$
with its weights $\wt(x,y,z,w)=(1,5,3,2)$.
\end{Claim}

\begin{proof}[Proof of Claim \textup{\ref{cl:lem:restate-excep-wbl}}]
First we remark that $x, z/x, w/x \in \cO_{X_1, Z_1}$
generate local coordinates of $X_1$ at $Z_1$,
that $y/x = -((z/x)^2+x^{N-2}(w/x)^N)$,
and that $F_3$ equals, as valuations,
the exceptional divisor obtained by the weighted blow-up of $X_1$
with its weights $\wt(x,z/x,w/x)=(1,2,1)$.
Thus it is trivial that
$(m_{\divis(x)}(F_3), m_{\divis(y)}(F_3),
m_{\divis(z)}(F_3), m_{\divis(w)}(F_3))
=(1,5,3,2)$.
Because any $v \in \cO_{X, P}$ has an expansion of a formal series
$v = v_1(x,z,w)+v_2(y,z,w)$, it is sufficient to prove that for any $i \ge 0$,
\begin{align*}
v = \sum_{(p,q,r,s) \in I_i} c_{pqrs}x^py^qz^rw^s 
\in g_{03*}\cO_{X_3}(-(i+1)F_3) \quad (c_{pqrs} \in \bC)
\end{align*}
implies $v=0$, where
\begin{align*}
I_i =\{(p,q,r,s) \in \bZ_{\ge 0}^4 | p+5q+3r+2s=i,\;\; p \textrm{ or } q=0\}.
\end{align*}
However, by replacing $v$ with $x^jv$ for a sufficiently large $j$,
we have only to show that for any $i \ge 0$,
\begin{align*}
v = \sum_{(p,q,r) \in J_i} c_{pqr}x^pz^qw^r
\in g_{03*}\cO_{X_3}(-(i+1)F_3) \quad (c_{pqr} \in \bC)
\end{align*}
implies $v=0$, where
$J_i =\{(p,q,r) \in \bZ_{\ge 0}^3 | p+3q+2r=i\}$.

Take any $v = \sum_{(p,q,r) \in J_i} c_{pqr}x^pz^qw^r
\in g_{03*}\cO_{X_3}(-(i+1)F_3)$.
Then $v = \sum_{(p,q,r) \in J_i} c_{pqr}x^{p+q+r}(z/x)^q(w/x)^r$.
Because $F_3$ equals, as valuations,
the exceptional divisor obtained by the weighted blow-up of $X_1$
with its weights $\wt(x,z/x,w/x)=(1,2,1)$,
it is enough to show that the weight of any monomial $x^{p+q+r}(z/x)^q(w/x)^r$
($(p,q,r) \in J_i$) with respect to
its weights $\wt(x,z/x,w/x)=(1,2,1)$ equals $i$.
But this is trivial by a direct calculation $(p+q+r)+2q+r=p+3q+2r=i$.
\end{proof}

Only the proof of $N=3$ remains.
Because of Lemma \ref{lem:excep} and properties of toric geometry,
we have only to show the following claim.

\begin{Claim}\label{cl:lem:restate-excep-exclude}
Consider an analytic germ of a c$A_1$ point
$o \in (xy+z^2+w^N=0) \subset \bC^4$ $(N \ge 4)$
and blow-up this with its weights $\wt(x,y,z,w)=(1,5,3,2)$.
Then the exceptional locus of this weighted blow-up is irreducible, and
the weighted blown-up analytic space is normal and
has a non-terminal singularity.
\end{Claim}

\begin{proof}[Proof of Claim \textup{\ref{cl:lem:restate-excep-exclude}}]
Direct calculation shows that its exceptional locus is isomorphic to
$(xy+z^2=0) \subset \bP(1,5,3,2)$ with weighted homogeneous coordinates
$x,y,z,w$, which is irreducible,
and that all singularities on the obtained analytic space
are one terminal quotient singularity of type $\frac{1}{5}(-1,3,2)$ and
one non-terminal singularity isomorphic to
$o \in (xy+z^2+w^{2N-6}=0) \subset \bC^4/\bZ_2(1,1,1,-1)$.
\end{proof}
\end{proof}

Now our problem is proving (\ref{lem:restate-excep}.1-2),
which will be shown in Lemmas \ref{lem:exclude}.1 and \ref{lem:system-locus}.
We show all the possible cases.

\begin{Lemma}\label{lem:a=4}
$a=4$, and the tower $X_n \to \cdots \to X_0$ is exactly one of the following.

\smallskip
\noindent \textup{(\ref{lem:a=4}.1)}
$(n,m,l)=(3,2,1)$, $N \ge 3$, $r=5$.\\
\textup{(\ref{lem:a=4}.2)}
$(n,m,l)=(4,2,2)$, $N \ge 5$, $r=5$.\\
\textup{(\ref{lem:a=4}.3)}
$(n,m,l)=(4,3,3)$, $N \ge 7$, $r=7$.
\end{Lemma}

\begin{proof}
Though $a=2$ or $4$ in case (\ref{th:fromKwk}.2.1.1),
$a=2$ is impossible because $n \ge 3$. Hence $a=4$.
By Corollary \ref{cor:last}.2, it is trivial that
the values of $n,m,l$ in (\ref{lem:a=4}.1-3) cover
all the possibility for $a=4$ and $1<m<n$.

Now we calculate the value of $r$ in each case
usuig Proposition \ref{prop:fromKwk}.3.
Because $a=4$ and $J=\{(r,2)\}$ $(r \ge 5)$,
Proposition \ref{prop:fromKwk}.3 implies that
\begin{align*}
\dim_{\bC}\fm_P/f_*\cO_Y(-3E) &= 3+\max\{0, 6-r\}, \\
\dim_{\bC}\fm_P/f_*\cO_Y(-4E) &= 4+\max\{0, 8-r\}.
\end{align*}
Thus we have only to the next claim.

\begin{Claim}\label{cl:lem:a=4}
\textup{(\ref{cl:lem:a=4}.1)}
In case \textup{(\ref{lem:a=4}.1)}, $\dim_{\bC}\fm_P/f_*\cO_Y(-3E) = 4$.\\
\textup{(\ref{cl:lem:a=4}.2)}
In case \textup{(\ref{lem:a=4}.2)}, $\dim_{\bC}\fm_P/f_*\cO_Y(-3E) = 4$.\\
\textup{(\ref{cl:lem:a=4}.3)}
In case \textup{(\ref{lem:a=4}.3)}, $\dim_{\bC}\fm_P/f_*\cO_Y(-4E) = 5$.
\end{Claim}

\begin{proof}[Proof of Claim \textup{\ref{cl:lem:a=4}}]
we will express $f_*\cO_Y(-iE)$'s in each case
under a suitable identification
$P \in X \cong o \in (xy+z^2+w^N=0) \subset \bC^4$.

\smallskip
\noindent (\ref{cl:lem:a=4}.1)
As in Claim \ref{cl:lem:restate-excep-id},
we may assume that $l_0 = F_1 \cap \divis(y)_{X_1} \cap \divis(z)_{X_1}$
and $Z_1 = l_0 \cap \divis(w)_{X_1}$. Then
\begin{align*}
f_*\cO_Y(-2E) &= (y,z,w) + \fm_P^2, \\
f_*\cO_Y(-3E) &= (v,y) + (z,w)\fm_P + \fm_P^3,
\end{align*}
where $v = t_zz + t_ww + t_{x^2}x^2$ for some $t_z, t_w, t_{x^2} \in \bC$
such that $t_z$ or $t_w$ is non-zero. This implies (\ref{cl:lem:a=4}.1).

\smallskip
\noindent (\ref{cl:lem:a=4}.2)
We may assume that $l_0 = F_2 \cap \divis(y)_{X_2} \cap \divis(z)_{X_2}$. Then
\begin{align*}
f_*\cO_Y(-2E) &= (x,y,z) + \fm_P^2, \\
f_*\cO_Y(-3E) &= (y,z) + (x)\fm_P + \fm_P^3.
\end{align*}
This implies (\ref{cl:lem:a=4}.2).

\smallskip
\noindent (\ref{cl:lem:a=4}.3)
We may assume that $l_0 = F_3 \cap \divis(y)_{X_3} \cap \divis(z)_{X_3}$. Then
\begin{align*}
f_*\cO_Y(-2E) &= (x,y,z) + \fm_P^2, \\
f_*\cO_Y(-3E) &= (x,y,z) + \fm_P^3, \\
f_*\cO_Y(-4E) &= (y,z) + (x)\fm_P + \fm_P^4.
\end{align*}
This implies (\ref{cl:lem:a=4}.3).
\end{proof}
\end{proof}

We exclude cases (\ref{lem:a=4}.2-3), which shows (\ref{lem:restate-excep}.1).
Moreover we determine the values of $c_i$'s in case (\ref{lem:a=4}.1).

\begin{Lemma}\label{lem:exclude}
\textup{(\ref{lem:exclude}.1)}
Neither case \textup{(\ref{lem:a=4}.2)} nor case
\textup{(\ref{lem:a=4}.3)} occurs.\\
\textup{(\ref{lem:exclude}.2)}
In case \textup{(\ref{lem:a=4}.1)}, $c_1=4/5$, $c_2=8/5$, $c_3=8/5$.
\end{Lemma}

\begin{proof}
We note that $m_E(F_i) \in \frac{1}{r}\bZ$ for any $i$.
Using Remark \ref{rem:not:discrep-mult}.2, for any $i$ we have
\begin{align*}
\alpha_{K_X}(F_i) - \sum_{1 \le j \le i} c_j = \alpha_{K_Y}(F_i)
= \alpha_{f^*K_X+4E}(F_i) = \alpha_{K_X}(F_i) - 4m_E(F_i).
\end{align*}
Hence $\sum_{1 \le j \le i} c_j = 4m_E(F_i) \in \frac{4}{r}\bZ$,
and thus $\forall c_i \in \frac{4}{r}\bZ$.

But on the other hand $c_i$'s satisfy the relations
in Remark \ref{rem:not:discrep-mult}.2 and Corollary \ref{cor:restr}.2.3.3.
Using them we know that there is no possibility for such $c_i$'s
in cases (\ref{lem:a=4}.2-3), and that
(\ref{lem:exclude}.2) is the only possibility in case (\ref{lem:a=4}.1).
\end{proof}

Now it is sufficient to deal with only case (\ref{lem:a=4}.1).
(\ref{lem:restate-excep}.2) comes from the following lemma,
and therefore we finish the proof of Proposition \ref{prop:main-excep}.
Let $l_0'$ be the strict transform of $l_0$ on $X_2$.

\begin{Lemma}\label{lem:system-locus}
\textup{(\ref{lem:system-locus}.1)}
Let $\sM_{F_1}$ be the linear system on $F_1 \cong Q_0$
obtained by the total pull-back of $\sM_{X_1}$
with the inclusion map $F_1 \hookrightarrow X_1$.
Then $\sM_{F_1}$ is a $0$-dimensional linear system consisting
of some multiple of $l_0$.\\
\textup{(\ref{lem:system-locus}.2)}
Let $\sM_{F_2}$ be the linear system on $F_2 \cong \bP^2$
obtained by the total pull-back of $\sM_{X_2}$
with the inclusion map $F_2 \hookrightarrow X_2$.
Then $\sM_{F_2}$ is a $0$-dimensional linear system consisting
of some multiple of $Z_2$.
\end{Lemma}

\begin{proof}
(\ref{lem:system-locus}.1)
Let $c$ be the multiplicity of $\sM_{F_1}$ along $l_0$,
and let $l \subset F_1$ be a general line
on $F_1 \cong Q_0 \subset \bP^3$. Then,
\begin{align*}
c/2= (cl_0 \cdot l)_{F_1}
\le (\mu \sM_{F_1} \cdot l)_{F_1} =-c_1(F_1 \cdot l)_{X_1}=4/5.
\end{align*}
On the other hand,
\begin{align*}
-c/2= (cl_0' \cdot l_0')_{{F_1}_{X_2}}
&\le (\mu \sM_{X_2} \cdot l_0')_{X_2} \\
&=-c_2(F_2 \cdot l_0')_{X_2}-c_1(F_1 \cdot l_0)_{X_1} \\
&= -c_2 + c_1 = -4/5.
\end{align*}
By these two inequalities, we obtain $c = 8/5$
and $(cl_0 \cdot l)_{F_1} = (\mu \sM_{F_1} \cdot l)_{F_1}$.
This shows (\ref{lem:system-locus}.1).

\smallskip
\noindent (\ref{lem:system-locus}.2)
Because Corollary \ref{cor:restr}.2.3.2 tells that $Z_2 \subset F_2$
is a line on $F_2 \cong \bP^2$,
we know that $g_3$ induces an isomorphism
${F_2}_{X_3} \cong F_2 \cong \bP^2$.
Let $\sM_{{F_2}_{X_3}}$ be the linear system on ${F_2}_{X_3} \cong \bP^2$
obtained by the total pull-back of $\sM_{X_3}$
with the inclusion map ${F_2}_{X_3} \hookrightarrow X_3$.
It is enough to prove that $\sM_{{F_2}_{X_3}} = \emptyset$.

Let $l \subset F_2$ be a general line on $F_2 \cong \bP^2$, and
let $l'$ be the strict transform of $l$ on $X_3$. Then
\begin{align*}
(\mu \sM_{{F_2}_{X_3}} \cdot l')_{{F_2}_{X_3}}
=-c_3(F_3 \cdot l')_{X_3}-c_2(F_2 \cdot l)_{X_2}
=-c_3+c_2=0,
\end{align*}
which shows that $\sM_{{F_2}_{X_3}} = \emptyset$.
\end{proof}


\section{General Case}\label{sec:gen}
In this section we treat the remaining general case (\ref{th:main-th}.1),
which corresponds to cases (\ref{th:fromKwk}.2.1.2) and (\ref{th:fromKwk}.2.2),
and our aim is the following, which terminates the proof of
Theorem \ref{th:main-th}.

\begin{Proposition}\label{prop:main-gen}
Assume that $f$ is of type \textup{(\ref{th:fromKwk}.2.1.2)}
or \textup{(\ref{th:fromKwk}.2.2)}.
Then $f$ is obtained by a suitable weighted blow-up
of type \textup{(\ref{th:main-th}.1)}.
\end{Proposition}

Throughout this section except Definition \ref{def:sp-surf} and
Proposition \ref{prop:sp-surf}, we assume that $f$ is
of type (\ref{th:fromKwk}.2.1.2) or (\ref{th:fromKwk}.2.2)
and struggle with Proposition \ref{prop:main-gen}.
We set $(r_1, r_2) = (1, r)$ in case (\ref{th:fromKwk}.2.2) in this section
because we want to treat both cases (\ref{th:fromKwk}.2.1.2) and
(\ref{th:fromKwk}.2.2) simultaneously.

First we restate the conclusion.

\begin{Lemma}\label{lem:restate-gen}
The following imply Proposition \textup{\ref{prop:main-gen}}.

\smallskip
\noindent \textup{(\ref{lem:restate-gen}.1)}
$l=m$.\\
\textup{(\ref{lem:restate-gen}.2)}
$N \ge 2a$.\\
\textup{(\ref{lem:restate-gen}.3)}
There exists an identification
$P \in X \cong o \in (xy+z^2+w^N=0) \subset \bC^4$
satisfying that $z \in f_*\cO_Y(-aE)$.
\end{Lemma}

\begin{proof}
We use the same idea as that in the proof of Lemma \ref{lem:restate-excep}.
First we remark that $a=n$ by (\ref{lem:restate-gen}.1)
and Corollary \ref{cor:last}.2.
By (\ref{lem:restate-gen}.3), we have an identification
$P \in X \cong o \in (xy+z^2+w^N=0) \subset \bC^4$
satisfying that $z \in f_*\cO_Y(-nE)$.
Moreover by (\ref{lem:restate-gen}.1) we may assume that
$Z_m = F_m \cap \divis(y)_{X_m} \cap \divis(z)_{X_m} \subset
F_m \cong Q_0 \subset \bP^3$.
We know that $(m_{\divis(x)}(F_n), m_{\divis(z)}(F_n), m_{\divis(w)}(F_n))
=(m,n,1)$.

\begin{Claim}\label{cl:lem:restate-gen-wbl}
Under the  above identification, $F_n$ equals, as valuations,
an exceptional divisor obtained by the weighted blow-up of $X$
with its weights $\wt(x,y,z,w)=(m,2n-m,n,1)$.
\end{Claim}

\begin{proof}[Proof of Claim \textup{\ref{cl:lem:restate-gen-wbl}}]
First we remark that $z/w^m, w \in \cO_{X_m, Z_m}$
generate local coordinates of $X_m$ at the generic point of $Z_m$,
that $x/w^m \in \cO_{X_m, Z_m}^\times$,
that $y/w^m = -(x/w^m)^{-1}((z/w^m)^2+w^{N-2m})$,
and that $F_n$ equals, as valuations,
the exceptional divisor dominating $Z_m$
which is obtained by the weighted blow-up of $X_m$
along $Z_m$ with its weights $\wt(z^m/w,w)=(n-m,1)$.
Thus we obtain
$(m_{\divis(x)}(F_n), m_{\divis(y)}(F_n),
m_{\divis(z)}(F_n), m_{\divis(w)}(F_n))
= (m,2n-m,n,1)$, considering (\ref{lem:restate-gen}.2).
Because any $v \in \cO_{X, P}$ has an expansion of a formal series
$v = v_1(x,z,w)+v_2(y,z,w)$, it is sufficient to prove that for any $i \ge 0$,
\begin{align*}
v = \sum_{(p,q,r,s) \in I_i} c_{pqrs}x^py^qz^rw^s 
\in g_{0n*}\cO_{X_n}(-(i+1)F_n) \quad (c_{pqrs} \in \bC)
\end{align*}
implies $v=0$, where
\begin{align*}
I_i =\{(p,q,r,s) \in \bZ_{\ge 0}^4 |
mp+(2n-m)q+nr+s=i,\;\; p \textrm{ or } q=0\}.
\end{align*}
However, by replacing $v$ with $x^jv$ for a sufficiently large $j$,
we have only to show that for any $i \ge 0$,
\begin{align*}
v = \sum_{(p,q,r) \in J_i} c_{pqr}x^pz^qw^r
\in g_{0n*}\cO_{X_n}(-(i+1)F_n) \quad (c_{pqr} \in \bC)
\end{align*}
implies $v=0$, where
$J_i =\{(p,q,r) \in \bZ_{\ge 0}^3 | mp+nq+r=i\}$.

Take any $v = \sum_{(p,q,r) \in J_i} c_{pqr}x^pz^qw^r
\in g_{0n*}\cO_{X_n}(-(i+1)F_n)$.
Then $v = \sum_{(p,q,r) \in J_i} c_{pqr}(x/w^m)^p(z/w^m)^qw^{mp+mq+r}$.
We note that $x/w^m \in \cO_{X_m,Z_m}^\times$.
Because $F_n$ equals, as valuations,
the exceptional divisor dominating $Z_m$ which is
obtained by the weighted blow-up of $X_m$ along $Z_m$
with its weights $\wt(z/w^m,w)=(n-m,1)$,
it is enough to show that the weight of any monomial
$(z/w^m)^qw^{mp+mq+r}$ ($(p,q,r) \in J_i$) with respect to
its weights $\wt(z/w^m,w)=(n-m,1)$ equals $i$.
But this is trivial by a direct calculation $(n-m)q+(mp+mq+r)=mp+nq+r=i$.
\end{proof}

There remains only proving that $m, n$ are coprime.
Because of Lemma \ref{lem:excep} and properties of toric geometry,
we have only to show the following claim.

\begin{Claim}\label{cl:lem:restate-gen-exclude}
Consider an analytic germ of a c$A_1$ point
$o \in (xy+z^2+w^N=0) \subset \bC^4$ $(N \ge 2n)$
and blow-up this with its weights $\wt(x,y,z,w)=(m,2n-m,n,1)$,
where $m, n$ are positive integers with $m<n$ and are not coprime. 
Then the exceptional locus of this weighted blow-up is irreducible, and
the weighted blown-up analytic space is normal and
has a non-terminal singularity.
\end{Claim}

\begin{proof}[Proof of Claim \textup{\ref{cl:lem:restate-gen-exclude}}]
Direct calculation shows that its exceptional locus is isomorphic to
$(xy+z^2=0)$ or $(xy+z^2+w^{2n}=0) \subset \bP(m,2n-m,n,1)$
with weighted homogeneous coordinates
$x,y,z,w$, which is irreducible,
and that all singularities on the obtained analytic space
are two non-terminal quotient singularities of types $\frac{1}{m}(-1,n,1)$
and $\frac{1}{2n-m}(-1,n,1)$,
and in the case $N \ge 2n+2$ furthermore
one terminal Gorenstein singularity isomorphic to
$o \in (xy+z^2+w^{N-2n}=0) \subset \bC^4$.
\end{proof}
\end{proof}

Our problem is proving (\ref{lem:restate-gen}.1-3).
For this we introduce one definition, which also makes sense
in more general situation as in Section \ref{sec:singRR}.

\begin{Definition}\label{def:sp-surf} 
An algebraic surface $P \in S \subset X$
is said to be \textit{special of type $s$},
where $s$ is a positive integer, if it satisfies the following conditions.

\smallskip
\noindent (\ref{def:sp-surf}.1)
$S$ is normal and has a Du Val singularity of type $A_s$ at $P$.\\
(\ref{def:sp-surf}.2)
$f^*S = S_Y + aE$.\\
\end{Definition}

A special surface has beautiful properties.

\begin{Proposition}\label{prop:sp-surf}
Let $P \in S \subset X$ be a special surface of type $s$,
and let $f_S$ be the induced morphism from $S_Y$ to $S$.
Then $S_Y$ is normal and $K_{S_Y} = f_S^*K_S$.
Especially, the minimal resolution of $S$ factors through $S_Y$.
\end{Proposition}

\begin{proof}
It is sufficient to show that $S_Y$ is normal and that $K_{S_Y} = f_S^*K_S$,
because these imply the last part of the statement.
We will prove them simultaneously.

Let $\nu \colon \widetilde{S_Y} \to S_Y$ be the normalization of $S_Y$.
First we calculate the dualizing sheaf $\omega_{S_Y}$ on $S_Y$.
Let $Y^o \subseteq Y$ be the Gorenstein locus of $Y$.
We remark that $Y \setminus Y^o$ is a finite set.
By the adjunction formula, we obtain that
\begin{align*}
\omega_{S_Y}|_{Y^o \cap S_Y}
&= \omega_Y(S_Y) \otimes_{\cO_Y} \cO_{S_Y}|_{Y^o \cap S_Y} \\
&= f_S^*(\omega_X(S) \otimes_{\cO_X} \cO_S)|_{Y^o \cap S_Y}
= f_S^*\omega_S|_{Y^o \cap S_Y}.
\end{align*}
On the other hand, we know that $\omega_{S_Y}$ is $(\textrm{S}_2)$,
that $S_Y \setminus (Y^o \cap S_Y) \subseteq S_Y$ is of codimension $\ge 2$,
and that $f_S^*\omega_S$ is invertible.
Thus we obtain $\omega_{S_Y}=f_S^*\omega_S$, and
our problem is reduced to only proving that $\nu$ is isomorphism.

Second we calculate the dualizing sheaf $\omega_{\widetilde{S_Y}}$
on $\widetilde{S_Y}$.
Grothendieck duality tells that
\begin{align*}
\omega_{\widetilde{S_Y}}&=
\cHom_{\cO_{S_Y}}(\nu_*\cO_{\widetilde{S_Y}}, \omega_{S_Y}) \\
&=\cHom_{\cO_{S_Y}}(\nu_*\cO_{\widetilde{S_Y}}, f_S^*\omega_S) \\
&=\cHom_{\cO_{S_Y}}(\nu_*\cO_{\widetilde{S_Y}}, \cO_{S_Y})
\otimes_{\cO_{\widetilde{S_Y}}} \nu^*f_S^*\omega_S,
\end{align*}
where the remark that $\omega_S$ is invertible induces the third equality.

Because $S$ is canonical,
the above equation shows that the conductor ideal sheaf
$\cHom_{\cO_{S_Y}}(\nu_*\cO_{\widetilde{S_Y}}, \cO_{S_Y})
\subseteq \cO_{\widetilde{S_Y}}$ has to equal $\cO_{\widetilde{S_Y}}$.
Hence $\nu$ is isomorphism.
\end{proof}

We come back to cases (\ref{th:fromKwk}.2.2) and (\ref{th:fromKwk}.2.1.2)
treated in this section.
In our situation, the type of any special surface must be higher.

\begin{Lemma}\label{lem:higher-type}
Let $P \in S \subset X$ be a special surface of type $s$.
Then $s \ge r_1+r_2-1$.
\end{Lemma}

\begin{proof}
First we give easy statements about a Du Val singularity of type $A_s$.

\begin{Claim}\label{cl:DuVal}
Let $P \in S$ be an algebraic germ \textup{(}resp.\@ an analytic germ\textup{)}
of a Du Val singularity of type $A_s$ $(s \ge 1)$,
let $f_S \colon (S_Y \supset E) \to (S \ni P)$
be a non-isomorphic partial resolution
factored through by the minimal resolution of $S$,
and let $P \in C \subset S$ be a general hyperplane section through $P$.\\
\noindent \textup{(\ref{cl:DuVal}.1)}
$C$ has its multiplicity $1$ along every prime component of $E$,
that is, $f_S^*C = C_{S_Y}+E$.\\
\textup{(\ref{cl:DuVal}.2)}
The set $C_{S_Y} \cap E$ consists of $2$ points, say $Q_1, Q_2$.
These $Q_1, Q_2$ are Du Val singularities of types $A_{s_1}, A_{s_2}$
with $s_1+s_2 < s$ $(s_1, s_2 \ge 0)$.
Here we define a Du Val singularity of type $A_0$ as a smooth point.\\
\textup{(\ref{cl:DuVal}.3)}
For $i=1, 2$,
the local intersection number $(C_{S_Y} \cdot E)_{S_Y, Q_i}$
equals $1/(s_i+1)$.
\end{Claim}

\begin{proof}[Proof of Claim \textup{\ref{cl:DuVal}}]
Let  $f \colon (T \supset F) \to (S \ni P)$ be the minimal resolution of $S$,
and let $g \colon T \to S_Y$ be the induced morphism.
$F=\sum_{1 \le i \le s}F_i$ is a chain of $(-2)$-curves $F_i$'s.
We order the indices $i$'s so that they are compatible with
the order of $F_i$'s in this chain.
It is fundamental to see that $f^*C = C_T + F$
and that $C_T$ intersects $F$ exactly at
a point, say $P_1$, on $F_1 \setminus F_2$ and at
a point, say $P_2$, on $F_s \setminus F_{s-1}$ transversally,
where we omit $\setminus F_2$ and $\setminus F_{s-1}$ if $s=1$.
Let $s_1$ (resp.\@ $s_2$) be the smallest non-negative integer such that
$F_{s_1+1}$ (resp.\@ $F_{s-s_2}$) is not contracted by $g$.
Then $Q_i=g(P_i)$ ($i=1, 2$) is a Du Val singularity of type $s_i$,
and $C_{S_Y} \cap E$ consists of $Q_1, Q_2$.
Because $g^*F_{s_1+1}=(s_1+1)^{-1}F_1+\others$
(resp.\@ $g^*F_{s-s_2}=(s_2+1)^{-1}F_s+\others$),
we have $(C_{S_Y} \cdot E)_{S_Y, Q_1} = 1/(s_1+1)$
(resp.\@ $(C_{S_Y} \cdot E)_{S_Y, Q_2} = 1/(s_2+1)$).
\end{proof}

We begin to prove Lemma \ref{lem:higher-type}.
We keep the notation $f_S \colon S_Y \to S$ in Proposition \ref{prop:sp-surf}.
Let $P \in H \subset X$ be a general hyperplane section on $X$ through $P$.
Then $P \in C=H|_S \subset S$ is also
a general hyperplane section on $S$ through $P$.
Because $m_P \neq f_*\cO_Y(-2E)$, we have $f^*H=H_Y+E$ and
$f_S^*C=H_Y|_{S_Y}+E|_{S_Y}$.
The support of $E|_{S_Y}$ is exactly the exceptional locus of $f_S$,
and $f_S$ is factored through by the minimal resolution of $S$
by Proposition \ref{prop:sp-surf}.
Thus by Claim \ref{cl:DuVal}.1, we obtain that $E|_{S_Y}$ is reduced and that
$H_Y|_{S_Y}=C_{S_Y}$, the strict transform of $C$ on $S_Y$.

We calculate the intersection number
of $C_{S_Y}$ and $E|_{S_Y}$ around $f_S^{-1}(P)$.
\begin{align*}
(C_{S_Y} \cdot E|_{S_Y})_{S_Y}
&=(H_Y \cdot E \cdot S_Y)_Y \\
&=((f^*H-E) \cdot E \cdot(f^*S-aE))_Y \\
&=aE^3=(1/r_1)+(1/r_2),
\end{align*}
where the last equality comes from Proposition \ref{prop:fromKwk}.2.

By Claim \ref{cl:DuVal}.2,
the set $C_{S_Y} \cap E|_{S_Y}$ consists of $2$ points, say $Q_1, Q_2$,
and thus
$(C_{S_Y} \cdot E|_{S_Y})_{S_Y,Q_1}+(C_{S_Y} \cdot E|_{S_Y})_{S_Y,Q_2}
=(1/r_1)+(1/r_2)$.
We may assume that
$(C_{S_Y} \cdot E|_{S_Y})_{S_Y,Q_1} \ge (C_{S_Y} \cdot E|_{S_Y})_{S_Y,Q_2}$.
Considering the set $I$ and Claim \ref{cl:DuVal}.3, we know that
$(C_{S_Y} \cdot E|_{S_Y})_{S_Y, Q_1}=1/r_1$ and
$(C_{S_Y} \cdot E|_{S_Y})_{S_Y, Q_2}=1/r_2$,
and that the local Gorenstein indices of $Q_1, Q_2$ are $r_1, r_2$.
Therefore by Claims \ref{cl:DuVal}.2-3,
we obtain that $Q_1, Q_2$ are Du Val singularities of
types $A_{r_1-1}, A_{r_2-1}$ with
$(r_1-1)+(r_2-1)<s$, that is, $r_1+r_2 \le s+1$.
\end{proof}

\begin{Remark}\label{rem:lem:higher-type}
The above proof tells that
$Y$ has exactly two non-Gorenstein singularities
in case \textup{(\ref{th:fromKwk}.2.1.2)}.
\end{Remark}

We obtain an upper-bound of the value of $a$.

\begin{Lemma}\label{lem:upper-bound}
$r_1+r_2 \ge 2a$.
\end{Lemma}

\begin{proof}
Proposition \ref{prop:fromKwk}.2 induces that $a(r_1r_2E^3)=r_1+r_2$.
Thus we have only to show that $a \neq r_1+r_2$
because of Proposition \ref{prop:fromKwk}.1.
$m_E(F_1) \in \bZ$ (resp.\@ $\frac{1}{r_1}\bZ$, $\frac{1}{r_2}\bZ$)
when the center of $F_1$ on $Y$ is not a non-Gorenstein point
(resp.\@ is the non-Gorenstein point of index $r_1$,
is the non-Gorenstein point of index $r_2$).
Like the proof of Lemma \ref{lem:exclude}, we obtain
$c_1 \in a\bZ$ (resp.\@ $\frac{a}{r_1}\bZ$, $\frac{a}{r_2}\bZ$). 
By this and Proposition \ref{prop:eval}.1 we have
$a$ (resp.\@ $\frac{a}{r_1}$, $\frac{a}{r_2}$) $<1$,
which implies that $a \neq r_1+r_2$.
\end{proof}

Combining Lemmas \ref{lem:higher-type} and \ref{lem:upper-bound},
we obtain a corollary.

\begin{Corollary}\label{cor:higher-type}
Let $P \in S \subset X$ be a special surface of type $s$.
Then $s \ge 2a-1$.
\end{Corollary}

Now we will prove (\ref{lem:restate-gen}.1-3)
by constructing special surfaces.

\begin{Lemma}\label{lem:id-gen}
There exists an identification
$P \in X \cong o \in (xy+z^2+w^N=0) \subset \bC^4$
satisfying that $m_{\divis(w)}(E) = 1$ and that $z+p(w) \in f_*\cO_Y(-aE)$
for some $p(w) \in \bigoplus_{i=1}^{a-1}\bC w^i \subset \bC[w]$.\\
\end{Lemma}

\begin{proof}
The following claim is inevitable.

\begin{Claim}\label{cl:lem:id-gen-eqn}
Let $N_i = \lfloor \frac{i}{r_1} \rfloor +1$, which is
the number of elements in the set
$I_i= \{(j,k) \in \bZ_{\ge 0}^2 | r_1j+k=i\}$.
Then for $0 \le i < a$,
\begin{align*}
\dim_{\bC}f_*\cO_Y(-iE)/f_*\cO_Y(-(i+1)E)=N_i.
\end{align*}
\end{Claim}

\begin{proof}[Proof of Claim \textup{\ref{cl:lem:id-gen-eqn}}]
We note that $r_2 \ge a$ by Lemma \ref{lem:upper-bound}.
Thus by calculation using Proposition \ref{prop:fromKwk}.3
as in \cite[Remark 2.8.1, (2.12), (2.15)]{Kwk00},
we have, for $1 \le i \le a$,
\begin{align*}
\dim_{\bC}\cO_X/f_*\cO_Y(-iE)=
i - \frac{1}{2} \min_{0 \le j < i}\{((1+j)r_1-2i)j\}.
\end{align*}
\cite[Lemma 2.9]{Kwk00} shows that the above dimension equals
$\sum_{s=0}^{i-1} N_s$, which implies the claim.
\end{proof}

We express $f_*\cO_Y(-iE)$'s explicitly using the above claim.

\begin{Claim}\label{cl:lem:id-gen-weak}
\textup{(\ref{cl:lem:id-gen-weak}.1)}
Take an identification
$P \in X \cong o \in (xy+z^2+w^N=0) \subset \bC^4$
satisfing that $m_{\divis(w)}(E) = 1$.
Then for $1 \le i \le \min\{r_1, a\}$,
\begin{align*}
f_*\cO_Y(-iE)=(x_i, y_i, z_i)+(w^i)
\end{align*}
for some $x_i=x+p_i^x(w), y_i=y+p_i^y(w), z_i=z+p_i^z(w)$
$(p_i^x(w), p_i^y(w), p_i^z(w) \in
\bigoplus_{j=1}^{i-1}\bC w^j \subset \bC[w])$.\\
\textup{(\ref{cl:lem:id-gen-weak}.2)}
Assume $r_1 < a$.\\
\textup{(\ref{cl:lem:id-gen-weak}.2.1)}
In \textup{(\ref{cl:lem:id-gen-weak}.1)},
$x_{r_1}$, $y_{r_1}$ or $x_{r_1}-y_{r_1}-2z_{r_1}
\not\in f_*\cO_Y(-(r_1+1)E) + (w^{r_1})$.\\
\textup{(\ref{cl:lem:id-gen-weak}.2.2)}
In \textup{(\ref{cl:lem:id-gen-weak}.1)},
assume that $x_{r_1} \not\in f_*\cO_Y(-(r_1+1)E) + (w^{r_1})$.
Under this situation, for $r_1 \le i \le a$,
\begin{align*}
f_*\cO_Y(-iE)=(y_i, z_i)+ \sum_{(j,k) \in \cup_{s \ge i} I_s}(x_{r_1}^j\!w^k)
\end{align*}
for some $y_i=y+p_i^y(x_{r_1}, w), z_i=z+p_i^z(x_{r_1}, w)$
$(p_i^y(x_{r_1}, w), p_i^z(x_{r_1}, w) \in \bigoplus_{s=1}^{i-1}
\bigoplus_{(j,k) \in I_s}\bC x_{r_1}^j\!w^k \subset \bC[x_{r_1}, w])$,
where $I_i$ is the set in Claim \textup{\ref{cl:lem:id-gen-eqn}}.
\end{Claim}

\begin{proof}[Proof of Claim \textup{\ref{cl:lem:id-gen-weak}}]
(\ref{cl:lem:id-gen-weak}.1)
We will construct $x_i, y_i, z_i$ inductively
starting with $x_1=x, y_1=y, z_1=z$.
Assume that we have constructed $x_i, y_i, z_i$ ($1 \le i < \min\{r_1, a\}$).
There exists a surjective map $\lambda_i$,
\begin{multline*}
\lambda_i \colon ((x_i, y_i, z_i)+(w^i))/(\fm_P(x_i, y_i, z_i)+(w^{i+1})) \\
\twoheadrightarrow f_*\cO_Y(-iE)/f_*\cO_Y(-(i+1)E).
\end{multline*}
By $i<\min\{r_1, a\}$ and Claim \ref{cl:lem:id-gen-eqn},
$\dim_{\bC}f_*\cO_Y(-iE)/f_*\cO_Y(-(i+1)E)=N_i=1$.
On the other hand because $m_{\divis(w)}(E) = 1$, we know that
$w^i$ generates $f_*\cO_Y(-iE)/f_*\cO_Y(-(i+1)E)$, and that
$x_i+t_xw^i, y_i+t_yw^i, z_i+t_zw^i \in \Ker \,\lambda_i$
for some $t_x, t_y, t_z \in \bC$.
Hence it is enough to put
$x_{i+1}=x_i+t_xw^i, y_{i+1}=y_i+t_yw^i, z_{i+1}=z_i+t_zw^i$.

\smallskip
\noindent (\ref{cl:lem:id-gen-weak}.2.1)
As in the above proof,
using $x_{r_1}, y_{r_1}, z_{r_1}$ in (\ref{cl:lem:id-gen-weak}.1)
we have a surjective map $\lambda_{r_1}$,
\begin{multline*}
\lambda_{r_1} \colon
((x_{r_1}, y_{r_1}, z_{r_1})+(w^{r_1}))
/(\fm_P(x_{r_1}, y_{r_1}, z_{r_1})+(w^{r_1+1})) \\
\twoheadrightarrow f_*\cO_Y(-r_1E)/f_*\cO_Y(-(r_1+1)E).
\end{multline*}
Dividing by $(w^{r_1})$, we have another surjective map
$\overline{\lambda}_{r_1}$,
\begin{multline*}
\overline{\lambda}_{r_1} \colon
(x_{r_1}, y_{r_1}, z_{r_1})/\fm_P(x_{r_1}, y_{r_1}, z_{r_1}) \\
\twoheadrightarrow f_*\cO_Y(-r_1E)/(f_*\cO_Y(-(r_1+1)E)+(w^{r_1})).
\end{multline*}
By Claim \ref{cl:lem:id-gen-eqn} and $m_{\divis(w)}(E) = 1$,
$\dim_{\bC}f_*\cO_Y(-r_1E)/(f_*\cO_Y(-(r_1+1)E+(w^{r_1}))=N_{r_1}-1=1$.
Hence $\dim_{\bC}\Ker \,\lambda_{r_1} = 3-1 = 2$,
which shows (\ref{cl:lem:id-gen-weak}.2.1).

\smallskip
\noindent (\ref{cl:lem:id-gen-weak}.2.2)
We will prove (\ref{cl:lem:id-gen-weak}.2.2) as in the proof of
(\ref{cl:lem:id-gen-weak}.1), constructing $y_i, z_i$ inductively
starting with $y_{r_1}, z_{r_1}$ in (\ref{cl:lem:id-gen-weak}.1).
Assume that we have constructed $y_i, z_i$ ($r_1 \le i < a$).
There exists a surjective map $\lambda_i$,
\begin{multline*}
\lambda_i \colon ((y_i, z_i)+
\sum_{(j,k) \in \cup_{s \ge i} I_s}(x_{r_1}^j\!w^k))
/(\fm_P(y_i, z_i)+ \sum_{(j,k) \in \cup_{s \ge i+1} I_s}(x_{r_1}^j\!w^k)) \\
\twoheadrightarrow f_*\cO_Y(-iE)/f_*\cO_Y(-(i+1)E).
\end{multline*}
We know that $x_{r_1}, w^{r_1}$ generate
$f_*\cO_Y(-r_1E)/f_*\cO_Y(-(r_1+1)E)$ because of the proof
of (\ref{cl:lem:id-gen-weak}.2.1).
Thus any non-zero element in
$\bigoplus_{(j,k) \in I_i}\bC x_{r_1}^j\!w^k \subset \bC[x_{r_1}, w]$,
which always decomposes into a product of
$w^{i - \lfloor \frac{i}{r_1} \rfloor r_1}$
and $\lfloor \frac{i}{r_1} \rfloor$ linear combinations of
$x_{r_1}, w^{r_1}$, has exactly its multiplicity $i$ along $E$.
This and Claim \ref{cl:lem:id-gen-eqn} imply that
$\{x_{r_1}^j\!w^k\}_{(j,k) \in I_i}$
generate $f_*\cO_Y(-iE)/f_*\cO_Y(-(i+1)E)$, and that
$y_i+t_i^y, z_i+t_i^z \in \Ker \,\lambda_i$
for some $t_i^y, t_i^z \in \bigoplus_{(j,k) \in I_i}\bC x_{r_1}^j\!w^k
\subset \bC[x_{r_1}, w]$.
Hence it is enough to put $y_{i+1}=y_i+t_i^y, z_{i+1}=z_i+t_i^z$.
\end{proof}

We will construct an identification in Lemma \ref{lem:id-gen}
using Claim \ref{cl:lem:id-gen-weak}.
It is easy that we can take an identification in (\ref{cl:lem:id-gen-weak}.1).
Lemma \ref{lem:id-gen} is trivial
if $a \le r_1$ by Claim \ref{cl:lem:id-gen-weak}.1.
If $r_1 < a$, by Claim \ref{cl:lem:id-gen-weak}.2.1 and an equation
$xy+z^2+w^N = (x-y-2z)y+(y+z)^2+w^N$, we may assume that
$x_{r_1} \not\in f_*\cO_Y(-(r_1+1)E)$ in the construction of
$x_{r_1}, y_{r_1}, z_{r_1}$ in (\ref{cl:lem:id-gen-weak}.1).
Then by Claim \ref{cl:lem:id-gen-weak}.2.2,
we obtain that $z_a = z+p_a^z(x+p_{r_1}^x(w), w) \in f_*\cO_Y(-aE)$.
We express $z_a$ as $z_a = z+p(w)+q(x,w)x$
($p(w) \in \bigoplus_{i=1}^{a-1}\bC w^i
\subset \bC[w]$, $q(x,w) \in \bC[x,w]$).
Thus it is sufficient to replace $y, z$ with
$y'=y-2q(x,w)z-q(x,w)^2x$, $z'=z+q(x,w)x$
because $xy+z^2+w^N=xy'+(z')^2+w^N$.
\end{proof}

Corollary \ref{cor:higher-type}, Lemma \ref{lem:id-gen}, and
the following lemma induce (\ref{lem:restate-gen}.1-3),
which terminates the proof of Proposition \ref{prop:main-gen} and therefore
also the proof of Theorem \ref{th:main-th} completely.

\begin{Lemma}\label{lem:constr-sp-surf}
\textup{(\ref{lem:constr-sp-surf}.1)}
Under the identification
$P \in X \cong o \in (xy+z^2+w^N=0) \subset \bC^4$ in
Lemma \textup{\ref{lem:id-gen}},
assume $N<2a$ or $p(w) \neq 0$.
Then there exists a special surface of type $s$ with $s<2a-1$.\\
\textup{(\ref{lem:constr-sp-surf}.2)}
Under the identification
$P \in X \cong o \in (xy+z^2+w^N=0) \subset \bC^4$ in
Lemma \textup{\ref{lem:id-gen}},
assume $N \ge 2a$, $p(w)=0$, and $l<m$.
Then there exists a special surface of type $2a-3$.
\end{Lemma}

\begin{proof}
(\ref{lem:constr-sp-surf}.1)
Take a surface $P \in S=\divis(z+p(w)+cw^a)$ for a general $c \in \bC$.
Then $P \in S \cong o \in (xy+(p(w)+cw^a)^2+w^N=0) \subset \bC^3$,
which is a Du Val singularity of type $A_s$,
where $s=\min\{2a, a+\ord \,p(w), \ord \,(p(w)^2+w^N)\}-1$.
Here $\ord \,q(w) = \inf\{i \in \bZ_{\ge0}|w^i \textrm{ divides } q(w)\}
\in \bZ_{\ge0} \cup \{+\infty\}$.
We remark that $s<2a-1$ if $N<2a$ or $p(w) \neq 0$.
Because $z+p(w) \in f_*\cO_Y(-aE)$ and
$m_{\divis(w)}(E)=1$, the multiplicity of $S$ along $E$ equals $a$.
Thus $P \in S \subset X$ is special of type $s$.

\smallskip
\noindent (\ref{lem:constr-sp-surf}.2)
We may assume that $l_0 = F_l \cap \divis(y)_{X_l} \cap \divis(z)_{X_l}$.
Since $l<m$,
$Z_l$ is a point on $l_0$ except the vertex point of $F_l \cong Q_0$.
Thus $Z_l = l_0 \cap \divis(tx+w^l)_{X_l}$ for some $t \in \bC$.
We note that $tx+w^l \in g_{0,l+1*}\cO_{X_{l+1}}(-(l+1)F_{l+1})
\subseteq f_*\cO_Y(-(l+1)E)$ because $Z_l \in \divis(tx+w^l)_{X_l}$.
Take a surface $P \in S=\divis(z+w^{a-l-1}(tx+w^l)+cw^a)$
for a general $c \in \bC$.
Then $P \in S \cong o \in (xy+(w^{a-l-1}(tx+w^l)+cw^a)^2+w^N=0)
\subset \bC^3$,
which is a Du Val singularity of type $A_{2a-3}$.
Because $z \in f_*\cO_Y(-aE)$, $tx+w^l \in f_*\cO_Y(-(l+1)E)$,
and $m_{\divis(w)}(E)=1$, the multiplicity of $S$ along $E$ equals $a$.
Thus $P \in S \subset X$ is special of type $2a-3$.
\end{proof}


\noindent
\small{
\textsc{Department of Mathematical Sciences, University of Tokyo, Komaba,
Meguro, Tokyo 153-8914, Japan}
}\\
\textit{E-mail address:}
\texttt{kawakita@ms.u-tokyo.ac.jp}
\end{document}